\documentclass [12pt] {article}
\usepackage{amsfonts}
\usepackage{amssymb}

\newcommand{\qed} {\hspace {0.1in} \rule {1.5mm} {3.5mm}}

\newtheorem{lemma}{Lemma}[section]

\newtheorem{theorem}{Theorem}
\newtheorem{remark}{Remark}[section]

\newtheorem{definition}{Definition}[section]

\def\ps{\Psi_n}

\def\G{\Gamma}
\def\map{\mathbb Map}

\def\<{\langle}
\def\>{\rangle}
\def\s{\sigma}
\def\proof{\smallskip\noindent{\bf Proof:} }

\def\to{\rightarrow}

\title{On sofic groups} 
\author{{\sc G\'abor Elek and Endre Szab\'o}
\cr Mathematical Institute of
the Hungarian Academy of Sciences\cr P.O. Box 127, H-1364 Budapest, Hungary\cr
elek@renyi.hu\cr
endre@renyi.hu}
\date{}
\begin{document}

\maketitle
\noindent{\bf Abstract.} Answering some queries of Weiss \cite{Wei},
we prove that the free product and amenable extensions of sofic groups
are sofic as well, and give an example of a finitely generated
sofic group that is not residually amenable. 
\vskip 0.2in
\newpage
\section{Introduction}
Sofic groups (originally: initially subamenable groups) were introduced
by Gromov \cite{Gro1}. They can be viewed as a common 
generalizations of amenable and residually finite groups. 
Our interest in sofic groups arouse when in \cite{ES}
we proved Kaplansky's direct finiteness conjecture for sofic groups.
First of all let us recall the definition of soficity.
For a finite set $A$ let $\map(A)$ denote the monoid of self-maps of
$A$ acting on the right. We use multiplicative notation for the action,
that is we write $a\cdot f$ for $f(a)$ and multiplication in $\map(A)$
works as usual:  
$a\cdot fg = (a\cdot f)\cdot g = g(f(a))$.
\noindent
Let $\epsilon\in (0,1)$ be a real number, then
we say that two elements $e,f\in\map(A)$ are {\em $\epsilon$-similar},
or $e\sim_\epsilon f$, 
if the number of points $a\in A$ with $a\cdot e\neq a\cdot f$ is at most
$\epsilon|A|$. We say that these $e,f$ are $(1-\epsilon)$-different, or
$e\not\sim_{1-\epsilon}f$, if they are not $(1-\epsilon)$-similar. 

\begin{definition}
\label{quasi-action}
Let $G$ be a group, $\epsilon\in(0,1)$ a real number and $F\subseteq G$ a
finite subset. 
An {\em $(F,\epsilon)$-quasi-action} of $G$ on a finite set $A$
is a function $\phi:G\to\map(A)$ with the following properties: 
\begin{enumerate}
\item[(a)]   For any two elements $e,f\in F$ the map $\phi(ef)$ is
             $\epsilon$-similar to $\phi(e)\phi(f)$.
\item[(b)]   $\phi(1)$ is $\epsilon$-similar to the identity map of $A$.
\item[(c)]   For each $e\in F\setminus\{1\}$ the map $\phi(e)$ is
        $(1-\epsilon)$-different from the identity map of $A$. 
\end{enumerate}
\end{definition}

\begin{definition}
\label{initial}
The group $G$ is {\em sofic} if for each number
$\epsilon\in(0,1)$ 
and any finite subset $F\subseteq G$ there exists an
$(F,\epsilon)$-quasi-action of $G$.
\end{definition}
\begin{remark} \
\label{quasi-action-remarks}
\begin{enumerate}
\item    If a countable group is sofic then there exists a
         countable sequence $(A_n,\phi_n)$ of quasi-actions
         such that for each finite
         set $F\subseteq G$ and each $\epsilon\in(0,1)$ the
         $(A_n,\phi_n)$ are $(F,\epsilon)$-quasi-actions
         for all sufficiently large index $n$. 
\item    It is enough to define the function $\phi$ of an
         $(F,\epsilon)$-quasi-action  on the subset
         $\{1\}\bigcup F\cdot F$, then one can choose an arbitrary
         extension to the whole of $G$.
\end{enumerate}
\end{remark}
\noindent
The obvious examples for sofic groups are the residually amenable 
groups. The goal of this paper is to answer some queries from the
survey of Weiss \cite{Wei}.
Namely, we prove that the class of sofic groups is closed
under free products and extensions by amenable groups and construct
an example of a finitely generated sofic group which is not residually
amenable. 
\section{The class of sofic groups}
\begin{lemma}
\label{good-action}
If a group $G$ is sofic then for each finite subset
$F\subseteq G$ and for each $\epsilon\in(0,1)$ there is an
$(F,\epsilon)$-quasi-action $\phi$ of $G$ on a finite set $A$
satisfying the following extra conditions:
\begin{enumerate}
\item[(b')]   $\phi(1)$ is the identity map of $A$, for each 
              $1\ne g\in G$ the map $\phi(g)$ is a 
              fixpoint free bijection
              and $\phi(g^{-1})=\phi(g)^{-1}$. 
\item[(c')]   For different elements $e,f\in F\bigcup\{1\}$ the map 
              $\phi(e)$ is $(1-\epsilon)$-different from $\phi(f)$.
\end{enumerate}
\end{lemma}

\proof
Let $F^{-1}$ denote the set of inverses of the elements of $F$,
and let
$\tilde F = \break
\bigg(F\bigcup F^{-1}\bigcup\{1\}\bigg)\cdot
\bigg(F\bigcup F^{-1}\bigcup\{1\}\bigg)$
denote the collection of products of all pairs from the set
$\bigg(F\bigcup F^{-1}\bigcup\{1\}\bigg)$.
\noindent
Let us choose an 
$(\tilde F,\epsilon/10)$-quasi-action
$(A,\phi)$ and then we define a new quasi-action $\psi$ on the set 
$A'=A\coprod A$ (the disjoint union of two copies of $A$) as follows:
\noindent
Since $':X\to X'=X\coprod X$ is a functor,
so $\map(A)$ has a natural action on $A'=A\coprod A$,
acting the same way on both copies of $A$.
Using this action we get from $\phi$ another 
$(\tilde F,\epsilon/10)$-quasi-action  
$\phi'=\phi\coprod\phi$ \ on the set $A'$.
We use this trick only to make sure that important subsets have even
number of elements.
\noindent
Let $\psi(1)$ be the identity map on $A'$.
For each element $g$ outside of $\tilde F$
we define $\psi(g)$ to be an arbitrary fixpoint free involution on $A'$.
Here we need that $A'$ has an even number of elements.  
This way condition (b') holds for these elements, and the other
conditions are certainly uneffected. 

For each element $e\in \tilde F\setminus\{1\}$
we define $\psi(e)$ and $\psi(e^{-1})$ together.
Let $A_e\subseteq A$ denote the difference of the fixpoint sets of
$\phi(e)\phi(e^{-1})$ and $\phi(e)$.
In other words, $A_e$ is the largest subset of $A$ 
such that $\phi(e)\phi(e^{-1})$ is the identity map on $A_e$,
but $\phi(e)$ has no fixpoint in it.
We define $A_{e^{-1}}$ similarily (with $e^{-1}$ in place of $e$), 
it is clear that
$A_{e^{-1}}=A_e\cdot\phi(e)$.

Now $\phi(e)$ and $\phi(e^{-1})$
are inverse bijections between $A_e$ and $A_{e^{-1}}$. 
Our new $\psi(e)$ will be equal to $\phi'(e)$ on $A_e'=A_e\coprod A_e$
and we extend it to the whole of $A'$ in two steps.
First we extend to 
$(A_{e^{-1}}\bigcup A_e)'$
via an arbitrary bijection 
$
(A_{e^{-1}}\setminus A_e)' \to (A_e\setminus A_{e^{-1}})'
$.
Then extend it further via any fixpoint free involution on
the complement $\big(A\setminus(A_{e^{-1}}\bigcup A_e)\big)'$.
(One can do it since the complement has an even number of elements.)
This $\psi(e)$ is a fixpoint free bijection.
Then we define $\psi(e^{-1})$ to be its inverse.
(Note that our construction is symmetric in $e$ and $e^{-1}$.)
In case $e$ is an involution we find that $\phi(e)$ is an involution
on $A_e$, hence our $\psi(e)$ is also an involution.
Hence this $\psi$ satisfies condition (b') and therefore also 
condition (b) of Definition~\ref{quasi-action}.

Let's pick an element $e\in \tilde F\setminus\{1\}$.
We know that $\phi(e^{-1})\phi(e)$ is
$\epsilon/10$-similar to $\phi(1)$, hence $\epsilon/5$-similar to the
identity map. The fixpoint set of $\phi(e)$ has at most
$\frac\epsilon{10}|A|$ elements.
So the size of $A_e$ is at least $(1-\frac3{10}\epsilon)|A|$,
therefore $\psi(e)$ and $\phi'(e)$ are $\frac3{10}\epsilon$-similar.
The same is true for $e=1$ as well.
\noindent
Now let $e,f$ be elements of $F\bigcup F^{-1}$. Then $\psi(e)$ resp.
$\psi(f)$ are $\frac3{10}\epsilon$-similar to $\phi'(e)$ resp. $\phi'(f)$.
Since $\psi(e)$ is a bijection,
we find that $\psi(e)\psi(f)$ is $\frac3{10}\epsilon$-similar to
$\psi(e)\phi'(f)$. Hence
$$
\psi(e)\psi(f)\sim_{\frac3{10}\epsilon}
\psi(e)\phi'(f)\sim_{\frac3{10}\epsilon}
\phi'(e)\phi'(f)\sim_{\epsilon/10}
\phi'(ef)\sim_{\frac3{10}\epsilon}
\psi(ef).
$$
Putting it all together we find that $\psi(e)\psi(f)$ is
$\frac7{10}\epsilon$-similar to $\phi'(ef)$
and $\epsilon$-similar to $\psi(ef)$. 
Thus $\psi$ satisfies condition (a) of Definition~\ref{quasi-action}.

Now let $e,f$ be different elements of $F$. We have shown above that
$\psi(e)\psi(f^{-1})$ is $\frac7{10}\epsilon$-similar to $\phi'(ef^{-1})$
and the latter is $(1-\epsilon/10)$-different from the identity map.
Hence $\psi(e)\psi(f^{-1})$ is $(1-\frac8{10}\epsilon)$-different from the
identity. Since $\psi(f^{-1})=\psi(f)^{-1}$, we see that $\psi(e)$ and
$\psi(f)$ are $(1-\frac8{10}\epsilon)$-different. This proves condition
(c'), and also condition (c) of Definition~\ref{quasi-action}.
In particular, $\psi$ is an $(F,\epsilon)$-quasi-action satisfying (b')
and (c').

\qed

\begin{theorem}\label{constructions}
The class of sofic groups is closed under the following
constructions:
\begin{enumerate}
\item  \label{sofic-easy-constructions} 
       direct products, subgroups, inverse limits, direct limits,
\item  \label{sofic-free-product}       
       free products,
\item  \label{sofic-extension}
       certain extensions: if $N\lhd G$, $N$ is sofic
       and $G/N$ is amenable then $G$ is also sofic.
\end{enumerate}
\end{theorem}

\proof
We start with the {\bf proof of \ref{sofic-easy-constructions}}.
Let $\{G_i\}_{i\in I}$ be sofic groups, let
$G=\prod_{i\in I}G_i$ and $F\subseteq G$ a finite subset, and fix a
number $\epsilon\in(0,1)$. Then there exists a finite subset $J\subseteq I$
such that the natural projection $\pi_J:G\to G_J=\prod_{j\in J}G_j$
is injective on the set $F\cup\{1\}$. Then each
$\big(\pi_J(F),\epsilon\big)$-quasi-action of $G_J$ is also an
$(F,\epsilon)$-quasi-action of $G$, hence it is enough to prove
soficity for finite direct products. So we assume that the index
set is $I=\{1,2,\dots n\}$. 
Let $F_i$ denote the image of the projection of $F$ into the
factor $G_i$, and choose some $(F_i,\epsilon)$-quasi-actions
$\phi_i:G_i\to\map(A_i)$.
We define the finite set $A=\prod_{i\in I}A_i$ and the quasi-action
$\phi:G\to\map(A)$ via the formula
$$
\Big(a_1,a_2,\dots a_n\Big)\cdot\phi(g) =
\Big(a_1\cdot\phi(g_1),a_2\cdot\phi(g_2),\dots a_n\cdot\phi(g_n)\Big)
$$
It is obviously an
 $(F,n\epsilon)$-quasi-action of $G$. This shows that
$G$ is a sofic group.

Quasi-actions of a group $G$ are also quasi-actions for its subgroups,
hence if $G$ is sofic then all subgroups are sofic as well.
Inverse limit of sofic groups is by definition a subgroup
of their product, and therefore it is also sofic.

Now let $\{G_i\}_{i\in I}$ be a directed system of sofic groups, let
$G=\lim_{i\in I}G_i$, let  $F\subseteq G$ be a finite subset, and fix a
number $\epsilon\in(0,1)$. Then there is an index $i\in I$ and a
finite subset $F_i\subseteq G_i$ such that the natural homomorphism
$\sigma_i:G_i\to G$ is a bijection $F_i\to F$. Let $G'$ denote the
image of $G_i$ in $G$ and choose a coset representating system
$s:G'\to G_i$ (i.e. it has the property that $\sigma_i(s(g))=g$ for
all $g\in G'$). Let $\phi:G_i\to\map(A)$ be an
$(F_i,\epsilon)$-quasi-action of $G_i$. Then we define
$\Phi:G\to\map(A)$ via the formula
$$
a\cdot\phi(g) = \left\{{{a\cdot\phi\big(s(g)\big)}\atop{a}}
   \kern 15pt {{\mbox{when }g\in G',\ a\in A\hfill}
               \atop{\mbox{when }g\in G\setminus G',\ a\in A\hfill}}
\right.
$$
It is clearly an $(F,\epsilon)$-quasi-action of $G$. Hence $G$ is also
sofic. We finished proving \ref{sofic-easy-constructions}.

\vskip 10pt\noindent
Now we turn to the {\bf proof of \ref{sofic-free-product}} of the theorem.

\begin{definition}
\label{shortest-decomposition}
Let $G$ and $H$ be groups. Then each element $g\neq1$ of the free
product $G*H$ has a unique {\em shortest decomposition} of the form
$g=g_1h_1\dots g_kh_k$ where $g_i\in G$, $h_i\in H$.
For $g=1$ we define the shortest decomposition to be $g=1\cdot1$.
On can easily see that a decomposition 
$g=g_1h_1\dots g_kh_k$ 
is {\em shortest} if and only if none of the factors are $1$ 
except possibly $g_1$ or $h_k$.
\end{definition}

It is enough to prove soficity for the free product of two groups. So let
$G,H$ be sofic groups. Let $F_G$ and $F_H$ be finite
subsets of $G$ and $H$ respectively, $\epsilon\in(0,1)$ and  fix a
positive integer $N$. 
We pick an $(F_G,\epsilon)$-quasi-action $\phi$ on a finite set $A$ and
an $(F_H,\epsilon)$-quasi-action $\psi$ on a finite set $B$, both
satisfying the conditions of Lemma~\ref{good-action}.
\noindent
Let $F\subseteq G*H$ be the subset consisting of all elements with
shortest decomposition of the form
$g_1h_1\dots g_kh_k$ where $g_i\in F_G$, $h_i\in F_H$ and $k\le N$.
Our goal is to construct an $(F,\epsilon)$-quasi-action of the
free product $G*H$. This is enough for proving \ref{sofic-free-product}.
First we define the {\em incidence graph} of two
partitions. Then we construct our quasi-action in two steps.

\begin{definition}
Let $\alpha$ and $\beta$ be partitions of a finite set $C$. 
The {\em incidence graph} of $\alpha$ and $\beta$ is a bipartite
graph, whose two sets of vertices consist of the classes of $\alpha$
and the classes of $\beta$, and the edges are the elements of $C$, 
each element $c\in C$ connects its
$\alpha$-class with its $\beta$-class.
\end{definition}

\vskip 5pt \noindent
{\bf Step 1.}
{\it We construct a set $C$ with two partitions $\alpha$ and $\beta$
     on it with the following properties:
     each $\alpha$-class has $|A|$ elements, each $\beta$-class has
     $|B|$ elements, an $\alpha$-class can meet a $\beta$-class in at
     most one element, and in the incidence graph of $\alpha$ and $\beta$
     each circle is longer than $2N$.
}
\vskip 5pt
First we choose a finite group $V$ generated by $A\times B$ such that
all relations among the generators are longer than $2N$. The
existence of such finite groups follows from the fact that
free groups are residually finite.
In formulas we shall use the notation $(\overline{a,b})$ for
generators of $V$.
Let $C=A\times B\times V$. Our $\alpha$-classes will be the subsets
of the form
$$
A[b,v]=
\left\{\Big(a,b,v\Big)\Big|a\in A\right\}
\kern 20pt
\hbox{for each }
b\in B,\ v\in V
$$
We define $\beta$ classes as the subsets of the form
$$
B[a,w]=
\left\{\Big(a,b,w\cdot(\overline{a,b})\Big)\Big|b\in B\right\}
\kern 20pt
\hbox{for each }
a\in A,\ w\in V
$$
Two classes $A_{b,v}$ and $B_{a,w}$ can meet only if 
$v=w\cdot(\overline{a,b})$, 
and in this case their intersection is the single element $(a,b,v)$.
Suppose now that the classes
$$
A[b_1,v_1],\ B[a_2,v_2],\ A[b_3,v_3],\ B[a_4,v_4],\dots
A[b_{2k-1},v_{2k-1}],\ B[a_{2k},v_{2k}]
$$
form a circle with minimal length in the incidence graph of $\alpha$
and $\beta$. Our goal is to show that $k>N$.
Assume indirectly that $k\le N$.
To simplify notation, we shall use indexes modulo $2k$,
hence $b_{2k+1}=b_1$ and $v_{2k+1}=v_1$.
The above criterion for meeting classes now reads:
$$
v_{2i-1}=v_{2i}\cdot(\overline{a_{2i},b_{2i-1}}),
\kern 10pt
v_{2i+1}=v_{2i}\cdot(\overline{a_{2i},b_{2i+1}})
\kern 20pt
\hbox{for }
1\le i\le k
$$
This means that $v_1, v_2,v_3\dots v_{2k+1}=v_1$ is a returning path
in the Cayley graph of $V$. Since this graph have no loops of length
$2k$, this path must return along itself. This implies in turn that
there is a turning point, i.e. $v_{j+2}=v_j$ for some index $j$.
If $j$ is odd then we find that
$$
v_j = 
v_{j+1}\cdot(\overline{a_{j+1},b_j}) =
v_{j+2}\cdot(\overline{a_{j+1},b_{j+2}})^{-1}
       \cdot(\overline{a_{j+1},b_j})
$$
$$
v_j = v_j\cdot(\overline{a_{j+1},b_{j+2}})^{-1}
         \cdot(\overline{a_{j+1},b_j})
$$
hence
$$
(\overline{a_{j+1},b_j}) = \overline{(a_{j+1},b_{j+2}})
\kern 20pt
\hbox{so}
\kern 20pt
b_j=b_{j+2}
$$
But then $A[b_j,v_j]=A[b_{j+2},v_{j+2}]$, the original circle
can be shortened: 
$$
A[b_1,v_1],B[a_2,v_2], \dots
A[b_j,v_j],B[a_{j+3},v_{j+3}],\dots B[a_{2k},v_{2k}]
$$ 
This contradicts the minimality of the length $2k$.
For even $j$ the same argument gives us 
$B[a_j,v_j]=B[a_{j+2},v_{j+2}]$,
which is again a contradiction.
In both cases we run into contradiction, hence $k\le N$ is impossible.
This proves our claims about $\alpha$, $\beta$ and $C$.

\vskip 10pt \noindent
{\bf Step 2.}
{\it Using the two partitions we construct an
     $(F,\epsilon)$-quasi-action of the free product $G*H$ on our $C$.
}
\vskip 5pt
First we build an $(F_G,\epsilon)$-quasi-action
$\phi':G\to\map(C)$.
Since the $\alpha$-classes have $|A|$ elements, 
we make an identification $C\approx A\times(C/\alpha)$
so that the $\alpha$-classes are the subsets of the form
$A\times\{t\}$. 
We define $\phi'$ to act on the first coordinate:
$$
(a,x)\cdot\phi'(g) = \Big(a\cdot\phi(g),x\Big)
\kern 20pt
g\in G,\kern 3pt 
(a,x)\in C\approx A\times(C/\alpha).
$$
Similarily we define the $(F_H,\epsilon)$-quasi-action
$\psi':H\to\map(C)$
via the identification $C\approx B\times(C/\beta)$.

Now let $g\in G*H$ be any element with shortest decomposition
$g=g_1h_1\dots g_kh_k$ (where $g_i\in G$ and $h_i\in H$ 
as in Definition~\ref{shortest-decomposition}).
We define
$$
\Phi(g_1h_1\dots g_kh_k) = 
\phi'(g_1)\psi'(h_1)\dots \phi'(g_k)\psi'(h_k)
$$
Clearly $\Phi(1)$ is the identity map on $C$.

Next we prove that if $1\ne g\in F$ then $\Phi(g)$ has no fixpoints. 
Our proof is indirect.
Assume that some $c_0\in C$ is a fixpoint of $\Phi(g)$. 
Then we define a sequence of elements of $C$ via induction:
$$
c_{2i+1} = c_{2i}\cdot\phi'(g_{i+1})
\hbox{ \ \ \ and \ \ \ }
c_{2i+2} = c_{2i+1}\cdot\psi'(h_{i+1})
$$
for all $0\leq i\leq k-1$.
Now let $A_i$ be the $\alpha$-class of $c_{2i}$ and $B_i$ denote the
$\beta$-class of $c_{2i+1}$. Since $\phi'(g_{i+1})$ respects the
$\alpha$-classes we see that $c_{2i+1}$ also belongs to $A_i$, hence
$c_{2i+1}$ connects $A_i$ and $B_i$ in the incidence graph of
$\alpha$ and $\beta$. Similarily, $c_{2i+2}$ connects $B_i$ with
$A_{i+1}$ in this incidence graph. Hence we found a path 
$A_0,B_0\dots A_{k-1},B_{k-1},A_k$ in the incidence graph.

It is clear from the definition of $\Phi$ that
$$
c_{2k} = c_0\cdot \phi'(g_1)\psi'(h_1)\dots \phi'(g_k)\psi'(h_k)=
c_0\cdot\Phi(g)=c_0,
$$ 
hence $A_{2k}=A_0$, our path returns to $A_0$. 
But $g\in F$, so $k\le N$. 
Since the graph cannot not contain such a short circle,
our path must turn back at some point.
Let's look at the first turning point:
either we find $A_i=A_{i+1}$ or we have $B_i=B_{i+1}$
for some index $0\le i\le k-2$, or we have $k=1$.
(Indeed, if $A_{k-1}=A_k$ is a turning point and $k>1$, 
then $A_1\dots A_{k-1}$ is a shorter returning path, 
so there is another turning point before.)
In the first case 
$$
\{c_{2i+1}\} =
A_i\bigcap B_{i} =
A_{i+1}\bigcap B_{i} =
\{c_{2i+2}\} =
\Big\{c_{2i+1}\cdot\psi'(h_{i+1})\Big\}.
$$
This is impossible since $\psi'(h_{i+1})$ have no fixpoint for
$0\le i\le k-2$.
Similarily, $B_i=B_{i+1}$ would imply that $c_{2i+2}=c_{2i+3}$
is a fixed point of $\phi'(g_{i+2})$, which is again impossible.
In the third case we have $k=1$, and $c_1=c_2=c_0$
is a fixed point of both $g_1$ and $h_1$, hence both
$g_1=h_1=1$. This is also impossible, because $g_1h_1 = g \neq 1$.
Either way we have got a contradiction,
so $\Phi(g)$ must be fixpoint free.

Now we verify that our $\Phi$ is an $(F,\epsilon)$-quasi-action.
$\Phi(1)$ is the identity map, hence condition (b) of 
Definition~\ref{quasi-action} is satisfied.
Moreover, if $g\in F$ is different from $1$ then
$\Phi(g)$ is fixpoint free, 
hence condition (c) holds as well. We need to show (a).
Let $g,g'\in F$ with their shortest decompositions 
$g=g_1h_1\cdots g_kh_k$ and $g'=g'_1h'_1\dots g'_lh'_l$.
We distinguish three cases according to $h_k$ and $g'_1$:

\vskip 5pt
{\bf First}, 
if $h_k\neq1$ and $g'_1\neq1$ then there is no
cancellation and the shortest decomposition of $gg'$ is 
$g_1h_1\cdots h_kg'_1\dots g'_lh'_l$.
Then by definition $\Phi(gg')=\Phi(g)\Phi(g')$.

\vskip 5pt
{\bf Second}, 
if $h_k=1$ and $g'_1=1$ then the shortest decomposition
of $gg'$ is shorter, 
it is $g_1h_1\cdots g_kh'_1\dots g'_lh'_l$. 
Then both $\psi'(h_k)$ and $\phi'(g'_1)$ are the identity map,
so again 
$\Phi(gg')=\Phi(g)\Phi(g')$.

\vskip 5pt
{\bf Third}, 
if only one of $h_k$ and $g'_1$ is $1$
then in the product of the shortest decompositions of $g$ and $g'$
one factor is 1 (so must be omitted), some cancellations may happen
(e.g. if $g'_1=1$ then $h_k$ and $h'_1$ are neighbors and may cancel
each other),
and finally we must collapse the last remaining factor of $g$ and the
first remaining factor of $g'$ into one 
(e.g. if $h_k=1$ and nothing cancels then $g_k$ and $g'_1$ become
neighbors, and both are from the same group $G$). 
The shortest decomposition of $gg'$ has the form
$$
gg'=g_1h_1\cdots h_{k-t-1}(g_{k-t}g'_{1+t})h'_{1+t}\dots g'_lh'_l
$$
or
$$
gg'=g_1h_1\cdots g_{k-t}(h_{k-t}h'_{1+t})g'_{2+t}\dots g'_lh'_l
$$
where $t$ is the number of cancellations, possibly 0.
When we multiply together the formulas for $\Phi(g)$ and $\Phi(g')$
the same cancellations will occur, so we must compare
$$
\Phi(g)\Phi(g')=\phi'(g_1)\cdots \psi'(h_{k-t-1})
    \underbrace{\phi'(g_{k-t})\phi'(g'_{1+t})}
           \psi'(h'_{1+t})\dots \psi(h'_l)
$$
with
$$
\Phi(gg')=\phi'(g_1)\cdots \psi'(h_{k-t-1})
   \overbrace{\phi'(g_{k-t}g'_{1+t})}
         \psi'(h'_{1+t})\dots \psi(h'_l)
$$
or we compare
$$
\Phi(gg')=\phi'(g_1)\cdots \phi'(g_{k-t})
   \underbrace{\psi'(h_{k-t}h'_{1+t})}
         \phi'(g'_{2+t})\dots \psi(h'_l)
$$
with
$$
\Phi(g)\Phi(g')=\phi'(g_1)\cdots \phi'(g_{k-t})
   \overbrace{\psi'(h_{k-t})\psi'(h'_{1+t})}
         \phi'(g'_{2+t})\dots \psi(h'_l)
$$
\vskip 5pt\noindent
In any case the factor $\phi'(g_{k-t}g'_{1+t})$ 
is $\epsilon$-similar to $\phi'(g_{k-t})\phi'(g'_{1+t})$,
the factor $\psi'(h_{k-t}h'_{1+t})$ 
is $\epsilon$-similar to $\psi'(h_{k-t})\psi'(h'_{1+t})$ ,
and all other factors $\phi'(g_i)$, $\psi'(h_i)$ are bijections,
hence $\Phi(gg')$ is $\epsilon$-similar to the product
$\Phi(g)\Phi(g')$.

Therefore condition (a) of Definition~\ref{quasi-action} holds, our
$\Phi$ is an $(F,\epsilon)$-quasi-action. This finish the proof of the
fact that the free products of sofic groups are sofic groups as well.

\vskip 10pt\noindent
Now we {\bf prove \ref{sofic-extension}} of the theorem.
Let $N\lhd G$ be groups such that $N$ is sofic and
$G/N$ is amenable, and choose a finite set $F\subseteq G$ and a number
$\epsilon\in(0,1)$. Our goal is to build an
$(F,3\epsilon)$-quasi-action of $G$. 

For elements $g\in G$ let $\overline g\in G/N$ denote their image.
Let $\sigma:G/N\to G$ be a section of the natural homomorphism 
$G\to G/N$, i.e. it has the property that $\overline{\sigma(h)}=h$ for
each $h\in G/N$, or equivalently, $g\sigma(\overline{g})^{-1}\in N$ for
all $g\in G$.
By Folner's theorem we can choose a nonempty finite subset 
$\overline{A}\subseteq G/N$ with the property
that $|\overline A\overline g\setminus\overline A|\le \epsilon|\overline A|$ for all
$g\in F$. Let $A=\sigma(\overline A)$.
Moreover, let $H=N\bigcap(A\cdot F\cdot A^{-1})$,
and choose an $(H,\epsilon)$-quasi-action
$\psi$ of the group $N$ on a finite set $B$.

We define the map $\Phi:G\to\map(B\times A)$ as follows:
$$
\Big(b\,,\,a\Big)\cdot\Phi(g) = \left\{{
        {\Big(b\cdot\psi\big(ag\;\sigma(\overline a\overline g)^{-1}\big)
          \,,\,\sigma(\overline a\overline g)\Big)}  
          \atop  {\Big(b,a\Big)}}
   \kern 15pt {{\mbox{when }\overline a\overline g\in\overline A
         \lower 11pt \hbox{}}
            \atop  {\mbox{otherwise}\raise 11pt \hbox{}\hfill}}
\right.
$$
It is well defined since $ag\,\sigma(\overline a\overline g)^{-1}\in N$. 
We shall prove that $\Phi$ is an $(F,3\epsilon)$-quasi-action.
Clearly $\psi\big(a\sigma(\overline a)^{-1}\big)=\psi(1)$ is
$\epsilon$-similar to the identity map of $B$, hence 
$\Phi(1)$ is $\epsilon$-similar to the identity map of $B\times A$, 
so (b) of Definition~\ref{quasi-action} holds. 

Next we check condition~(c).
Let $e\in F\setminus\{1\}$, and assume first that $e\notin N$.
Then $\overline e\neq 1$, hence $a\ne\sigma(\overline a\overline e)$ 
for all $a\in G$. Therefore 
$$
\Big(b,a\Big)\cdot\Phi(e) =
\Big(b\cdot\psi(\dots),\sigma(\overline a \overline e) \Big)\ne
\Big(b,a\Big)
\kern 20pt
\hbox{ \ whenever \ }
\overline a\overline e\in\overline A
$$
Hence $(b,a)$ is a fixpoint of $\Phi(e)$ only if 
$\overline a\overline e\notin\overline A$. 
By assumption the number of such $\overline a$ is at most
$\epsilon|A|$, and each $\overline a$ gives $|B|$
fixpoints. So the total number of fixpoints is at most 
$\epsilon|A|\cdot|B|$.
Hence in this case $\Phi(e)$ is $(1-\epsilon)$-different from the
identity. 

Assume next that $e\in N\bigcap F\setminus\{1\}$.
In this case $\overline e=1$ hence
$$
\Big(b,a\Big)\cdot\Phi(e) =
\Big(b\cdot\psi(aea^{-1}),a\Big)
\kern 20pt
\hbox{for all }
b\in B,\ a\in A
$$
Hence $(b,a)$ is a fixpoint of $\Phi(e)$ if and only if $b$ is a
fixpoint of $\psi(aea^{-1})$. 
But $aea^{-1}\in H\setminus\{1\}$, so for a fixed $a\in A$ there are at most
$\epsilon|B|$ such fixpoints, hence altogether there are at most 
$\epsilon|A|\cdot|B|$ fixpoints, and thus 
$\Phi(e)$ is again $\epsilon$-different from the identity map.
This proves condition (c) of Definition~\ref{quasi-action}.

Finally let $e,f\in F$ and suppose that for certain $a\in A$ the
elements $\overline a\overline e$ and $\overline a\overline e\overline f$ are both in $\overline A$.
This assumption holds for all but $2\epsilon|A|$ values of $a$.
Then
$$
\Big(b,a\Big)\cdot\Phi(e)\Phi(f) =
\Big(b\cdot\psi\big(ae\;\sigma(\overline a\overline e)^{-1}\big)
          \,,\,\sigma(\overline a\overline e)\Big)\cdot\Phi(f) =
$$
$$
=
\Big(b\cdot\psi\big(ae\;\sigma(\overline a\overline e)^{-1}\big)
      \cdot\psi\big(\sigma(\overline a\overline e)f\;\sigma(\overline a\overline e\overline f)^{-1}\big)
          \,,\,\sigma(\overline a\overline e\overline f)\Big)
$$
and
$$
\Big(b,a\Big)\cdot\Phi(ef) =
\Big(b\cdot\psi\big(aef\;\sigma(\overline a\overline e\overline f)^{-1}\big)
          \,,\,\sigma(\overline a\overline e\overline f)\Big)
$$
But our assumption imply that 
$\sigma(\overline a\overline e),\sigma(\overline a\overline e\overline f)\in A$,
hence $ae\sigma(\overline a\overline e)^{-1}$ and 
$\sigma(\overline a\overline e)f\sigma(\overline a\overline e\overline f)^{-1}$ are elements
of $H$. Since $\psi$ is an $(H,\epsilon)$-quasi-action, we find that
with the exception of at most $\epsilon|B|$ values of $b$
$$
b\cdot\psi\big(ae\sigma(\overline a\overline e)^{-1}\big)
      \cdot\psi\big(\sigma(\overline a\overline e)f\sigma(\overline a\overline e\overline
      f)^{-1}\big) =
b\cdot\psi\big(aef\sigma(\overline a\overline e\overline f)^{-1}\big)
$$
Putting it together,
we have 
$$
(b,a)\cdot\Phi(e)\Phi(f) =
(b,a)\cdot\Phi(ef)
\kern 20pt 
\mbox{except for }
\left\{{
{2\epsilon|A|\mbox{ values of }a\hfill}
\atop { \epsilon|B|\mbox{ values of }b\hfill}}
\right.
$$
Therefore $\Phi(e)\Phi(f)$ is $3\epsilon$-similar to $\Phi(ef)$,
condition (a) of Definition~\ref{quasi-action} holds.

We proved all three conditions of Definition~\ref{quasi-action},
hence $\Phi$ is indeed an $(F,3\epsilon)$-quasi-action.
Hence \ref{sofic-extension} is proven.
This completes the proof of the theorem.

\qed

\begin{theorem}
\label{locally-residually-amenable}
Free product of locally residually amenable groups are sofic. 
In particular, free product of residually finite groups and amenable
groups are sofic.
\end{theorem}

\proof.
The one element group is sofic. Hence residually finite groups
are sofic by (c) of Theorem~\ref{constructions}.
Subgroups of direct products of amenable groups are sofic by
(a) of Theorem~\ref{constructions}. Residually amenable groups are
subgroups of the product of their amenable factors, therefore they
are sofic. Locally residually amenable groups are the direct
limit of residually amenable groups, hence they are also sofic
. Finally, free product of locally residually amenable
groups are sofic by (b) of Theorem~\ref{constructions}.

\qed
\section{A finitely generated non-residually amenable sofic group}
In \cite{VG}, the authors constructed a finitely generated amenable
LEF-group that was not residually finite. Modifying their construction
we show the existence of a finitely generated sofic group that is
not residually amenable.
Let $K$ be a hyperbolic, residually finite group with Kazhdan's Property (T)
\cite{Gro2}. Let $P$ be the set of all permutations of $K$ which move
only finitely many elements. Consider the group $Q$ generated by $P$
and the left translations by the elements of $K$. Note that $Q$ is the
semidirect product of $K$ and the locally finite group $P$.
\begin{theorem}
$Q$ is a finitely generated sofic group, but it is not residually amenable.
\end{theorem}
\proof Let $S$ be symmetric generating system for $K$. Denote by $T_S$ the
set of transpositions in the form $(1,s)$, where $s\in S$.
\begin{lemma}
\label{ujl1}
The set $W=S\cup T_S$ generate the group $Q$. \end{lemma}
\proof
Consider the right Cayley graph $\G$ of $K$ with respect to the generating
system $S$. Note that any transposition in the form $(g,gs)$, $g\in K$,
$s\in S$ can be written as $g\cdot(1,s)\cdot g^{-1}$. If $g,h\in K$ then
pick a shortest path 
$$ g\to gs_{i_1}\to gs_{i_1}s_{i_2}\to\dots \to gs_{i_1}s_{i_2}\dots s_{i_n}
=h$$
in the Cayley graph.
Any transposition in the form $(gs_{i_1}s_{i_2}\dots s_{i_k},
gs_{i_1}s_{i_2}\dots s_{i_k}s_{i_{k+1}})$ is generated by $W$ by our
previous observation, hence the transposition $(g,h)$ is generated by $W$
as well. This finishes the proof of our lemma. \qed
\begin{lemma}
\label{ujl2}
The group $Q$ is sofic.
\end{lemma}
\proof
Let $F_n\subseteq Q$ be the finite set of elements in the form of $k\s$,
where $k\in B_n(K)\subseteq\G$ and $\s\in P$ moves only the elements of
$B_n(K)$. Clearly, $F_n\cdot F_n\subseteq F_{2n}$ and $\cup_{n=1}^\infty F_n
=Q$, hence in order to prove that $Q$ is sofic, it is enough to construct
for any $n$ an injective map $\Psi_n$ from $F_{2n}$ to a finite group $H_n$
such that if $f,g\in F_n$, then $\Psi_n(f)\Psi_n(g)=\Psi_n(fg)$.

\noindent
Let $N_n$ be a normal subgroup of $K$ of finite index such that
$N_n\cap B_{10n}(K)=\{1\}\,.$
Such normal subgroup must exist since $K$ is residually finite.
Let $\tau:K\to K/{N_n}$ be the quotient homomorphism. The map
$\Psi_n:F_{2n}\to H_n= Sym(N_n)$ is defined as follows:

\noindent
$\Psi_n(k\sigma)=\Psi_n(k)\cdot\Psi_n(\s)$, where
$\ps(k)$ is the left translation by $\tau(k)$ and
$\ps(\s)(\tau(a))=\tau(\s(a))$ if $a\in B_{2n}(K)$, otherwise
$\ps(\s)(b)=b$. Obviously, $\ps$ is injective.

\noindent
Now let $s\in B_{5n}(K)$, $f=k_1\sigma_1\in F_1, g= k_2\s_2\in F_1$.
Then,
$$\ps(f)\ps(g)(\tau(s))=\ps(k_1\s_1)\ps(k_2)\tau(\s_2(s))=$$
$$=\ps(k_1)\tau(\s_1(k_2(\s_2(s)))=
\tau(k_1 \s_1  k_2 \s_2 (s))\,.$$
On the other hand,
$$\ps(k_1\s_1 k_2\s_2) (\tau(s))=
\ps(k_1k_2(k_2^{-1}\s_1  k_2)\s_2(\tau(s))=
\tau(k_1k_2) \tau((k_2^{-1} \s_1  k_2) \s_2)(s)))=$$
$$=\tau(k_1\s_1  k_2  \s_2 (s)\,.$$
Hence $\ps(f)\ps(g)$ coincides with $\ps(fg)$ on $\tau(B_{5n}(K))$.
However, if $l\notin \tau(B_{5n}(K))$ then
$$\ps(k_1\s_1)\ps(k_2\s_2)(l)=\tau(k_1 k_2)(l)$$
and
$$\ps(k_1\s_1 k_2\s_2)(l)=\tau(k_1 k_2)(l)\,.$$
Hence $\ps(f)\ps(g)=\ps(fg)$.\qed
\begin{lemma}
The group $Q$ is not residually amenable.
\end{lemma}
\proof
Let $A$ denote the simple subgroup of even permutations in $P$. 
If $Q$ were residually amenable then there exists a homomorphism
$\phi:Q\to M$, where $M$ is amenable, $\phi(t)\neq 1, \phi(a)\neq 1$,
$t$ is a non-torsion element in the
subgroup $K$ and $a$ is an even permutation. Note that such $t$ must
exist since $K$ is hyperbolic.
By the Kazhdan's property (T) the image of $K$ must be finite, since
any amenable quotient of a Kazhdan group is finite.
On the other hand, $\phi$ must be injective on $A$ since $A$ is simple
and $\phi(a)\neq 1$. Clearly, $t^{-n}a t^{-n}\neq a$ for any $n$.
However, if $n$ is the rank of $\phi(K)$, $\phi(t^n)=1$. Therefore,
$\phi(t^{-n}a t^{-n})=\phi(a)$ in contradiction with the
injectivity of $\phi$ on the subgroup $A$.
This finishes the proof of our lemma and of Theorem 3. as well. \qed

\end{document}